\documentclass[12pt]{amsart}

\usepackage{graphicx}
\usepackage{amsmath}
\usepackage{amsfonts}
\usepackage{amssymb}
\usepackage{latexsym}
\newtheorem{theorem}{Theorem}[section]                      
                       
\newtheorem{corollary}{Corollary}[section]  

\newtheorem{remark}{Remark}[section]

\theoremstyle{remark}

\input epsf

\begin{document}

\title[Lower Schwarz-Pick estimates and angular derivatives]
{Lower Schwarz-Pick estimates and angular derivatives}

\author[J. Milne Anderson and Alexander Vasil'ev]{J. Milne Anderson$^{1)}$ and Alexander Vasil'ev$^{2)}$}
\thanks{$^{1)}$Supported by the Leverhulme Trust (U.~K.)}
\thanks{$^{2)}$Supported by 
the grant of the University of Bergen (Norway)}
\address{J.~M.~Anderson: Department of Mathematics, University College London, Gower Street, London WC1E~6BT, U.K.}
\address{A.Vasil'ev: Department of Mathematics, University of Bergen, Johannes Brunsgate 12, Bergen 5008, Norway}
\email{alexander.vasiliev@uib.no}

\begin{abstract}  
The well-known Schwarz-Pick lemma states that any analytic mapping $\phi$ of the unit disk $U$ into itself satisfies the inequality $$|\phi'(z)|\leq \frac{1-|\phi(z)|^2}{1-|z|^2}, \quad z\in U.$$ This estimate remains the same if we restrict ourselves to univalent mappings. The lower estimate is $|\phi'(z)|\geq 0$ generally or  $|\phi'(z)|> 0$ for univalent functions. To make the lower estimate non-trivial we consider univalent functions and fix the angular limit and the angular derivative at some points of the unit circle. In order to obtain sharp estimates we make use of the reduced modulus of a digon.
\end{abstract}

\maketitle

\section{Introduction}

Let $U$ denote the unit disk; $U=\{z:\,\,|z|<1\}$. Pick's invariant form of the Schwarz lemma states that $$|\phi'(z)|\leq \frac{1-|\phi(z)|^2}{1-|z|^2},$$ for any analytic mapping $\phi:U\to U$ and $z\in U$. This estimate is obviously sharp
and remains true when we restrict ourselves to univalent $\phi$. Higher derivatives $|\phi^{(n)}|$ were estimated as well, see \cite{AndRovnyak, Ruscheweyh}. The lower estimates for $|\phi'(z)|$ are trivial. Two different reasons may be given to explain this phenomenon. The first one is the possible branching of an analytic mapping $\phi^{-1}$ at the point $\phi(z)$. The second is more appropriate to univalent mappings. One may choose a sequence of univalent mappings $\phi_n$ collapsing into the point $z$ as $n\to\infty$, and therefore, the conformal radius $|\phi'_{n}(z)$ tends to 0. To avoid such  behaviour one may prescribe a certain value to $\phi$ at another point of $U$ or at the boundary $\partial U$ in the angular sense. In the first case
this leads to bounded Montel's functions. However, even this is not sufficient to overcome collapsing. So we propose to fix also the derivative or the angular derivative at this extra point. The case of a fixed internal point was discussed in \cite{Vas02} in a slightly different context. However, letting the internal point tend to the boundary $\partial U$ gives no result and the boundary case must be treated separately. In this note we consider the class of conformal homeomorphisms $\phi$ of $U$
such that the angular limit $\angle \lim\limits_{z\to 1}\phi(z)=1$ exists and such that the angular limit $\angle \lim\limits_{z\to 1}\phi'(z)$ exists and equals $\beta$, $|\beta|<\infty$. In this case the mapping $\phi(z)$ is conformal at the boundary point 1.
Under these restrictions we obtain sharp lower estimates of $|\phi'(z)|$, $z\in U$. This result is a simple modification of \cite{Solynin} or \cite[Theorem 3.5.1]{VasBOOK} but we prove it in Section 4 to illustrate the method. The case of $n$ fixed boundary points is dealt with in Section 5. In both cases the estimates are sharp, with the extremal functions being given as the solution of a differential equation.

We thank Pietro  Poggi-Corradini and  Robert Burckel who informed us about the paper \cite{Unkelbach} where the simplest case of Theorems \ref{theoremA} ($\phi(0)=0$) and \ref{theorem} ($n=1$) was proved.

\section{Angular limits and angular derivatives}

In this section we present a short overview of the properties of angular limits and angular derivatives. For more information we refer to \cite{Pomm1, Pomm2}.
We say that a conformal mapping $\phi:\,U\to \mathbb C$ has the angular limit $\angle \lim\limits_{z\to \zeta}\phi(z)=\alpha$ at a point $\zeta\in \partial U$,
if $\lim\limits_{z\to \zeta,\,\,z\in \Delta_{\zeta}}\phi(z)=\alpha$ for any Stolz angle $\Delta_{\zeta}$ centered on $\zeta$.
If the limit $\lim\limits_{z\to \zeta}\phi(z)=\alpha$ exists for all $z\in U$, then $\phi$ becomes continuous at $\zeta$ as a function in $U\cup\{\zeta\}$. The angular limit $\alpha$ is a principal point of the image of the prime end of the mapping $\phi$ at the point $\zeta$. It is known that the angular limit $f(\zeta)=\alpha$ exists for almost all $\zeta \in \partial U$, moreover, the exceptional set in $\partial U$ has capacity zero. But $\phi$ is continuous at $\zeta$ only in some restricted cases. We say that $\phi$ has the angular derivative  $\phi'(\zeta)=\beta$ at the point $\zeta\in \partial U$ if the finite angular limit $f(\zeta)=\alpha$ exists and if
$$
\angle\lim\limits_{z\to\zeta}\frac{f(z)-\alpha}{z-\zeta}=\beta.
$$
The angular derivative exists if and only if the analytic function $\phi'(z)$ has the angular limit $\angle\lim\limits_{z\to\zeta}\phi'(z)=\beta$. Generally, very little may be said about the existence
of the angular derivatives. However, the Julia-Wolff theory implies that in the case of $\phi:\,U\to U$, the angular
derivative $\phi'(\zeta)$ exists (but perhaps may be infinite) at all points $\zeta\in \partial U$ where the angular limit $\phi(\zeta)$ exists and $|\phi(\zeta)|=1$. Furthermore, the mapping at the point $\zeta$ may be conformal ($0<|\phi'(\zeta)|<\infty$) or twisting. The McMillan Twist Theorem states that $\phi$ is conformal for almost all such points.  By the results of Julia, Carath\'eodory, and Wolff (see, e.g., \cite[page 57]{Nevanlinna} and \cite[page 306]{Pomm75}), if $\zeta$ is a boundary fixed point $\phi(\zeta)=\zeta$, then $\phi'(\zeta)$ is real.

\section{Moduli and extremal partitions}

\subsection{The reduced modulus of a digon}

Let $D$ be a hyperbolic simply connected domain in $\mathbb C$
with two finite fixed boundary points  $a$, $b$ (maybe with the
same support) on its piecewise smooth boundary. It is called a
digon. Denote by $S(a,\varepsilon)$ a region that is the
connected component of $D\cap\{|z-a|<\varepsilon\}$ with the point
$a$ in its border. Denote by $D_{\varepsilon}$ the domain
$D\setminus \{S(a,\varepsilon_1)\cup S(b,\varepsilon_2)\}$ for
sufficiently small $\varepsilon_{1,2}$ such that there is a curve in $D_{\varepsilon}$
connecting the opposite sides on $S(a,\varepsilon_1)$ and
$S(b,\varepsilon_2)$. Let $M(D_{\varepsilon})$ be the modulus of
the family of paths in $D_{\varepsilon}$ that connect the boundary
arcs of $S(a,\varepsilon_1)$ and $S(b,\varepsilon_2)$ when lie in
the circumferences $|z-a|=\varepsilon_1$ and $|z-b|=\varepsilon_2$
(we choose a single arc in each circle so that both arcs can be
connected in $D_{\varepsilon}$). If the limit
\begin{equation}
m(D,a,b)=\lim\limits_{\varepsilon_{1,2}\to 0}\left(
\frac{1}{M(D_{\varepsilon})}+
\frac{1}{\varphi_a}\log\,\varepsilon_1+\frac{1}{\varphi_b}\log\,\varepsilon_2\right),\label{eq:moddig}
\end{equation}
exists, where $\varphi_a=\sup\,\Delta_a$ and
$\varphi_b=\sup\Delta_b$ are the inner angles and $\Delta_{a,b}$
is the Stolz angle inscribed in $D$ at $a$ or $b$ respectively,
then $m(D,a,b$ is called the  reduced modulus of the digon $D$.
Various conditions guarantee the existence of this modulus,
whereas even in the case of a piecewise analytic boundary there are
examples \cite{Sol7} which show that this is not always the case.
The existence of the limit (\ref{eq:moddig}) is a local
characteristic of the domain $D$ (see \cite{Sol7}, Theorem 1.2).
If the domain $D$ is conformal (see the definition in \cite[page 80]{Pomm1}) at the points $a$ and $b$, then (\cite{Sol7}, Theorem
1.3)  the limit (\ref{eq:moddig}) exists. More generally, suppose
that there is a conformal map $f(z)$ of the domain
$S(a,\varepsilon_1)\subset D$ onto a circular sector, so that the
angular limit $f(a)$ exists which is thought of as a vertex of
this sector of angle $\varphi_a$. If the function $f$ has the
finite non-zero angular derivative $f'(a)$ we say that the domain
$D$ is also conformal at the point $a$ (compare \cite[page 80]{Pomm1}). If the digon $D$ is conformal at the points $a$, $b$
then  the limit (\ref{eq:moddig}) exists (\cite{Sol7}, Theorem
1.3). It is noteworthy that Jenkins and Oikawa \cite{JenOik}
in 1977 applied extremal length techniques to study the behaviour
of a regular univalent map at a boundary point. Necessary and
sufficient conditions were given for the existence of a finite
non-zero angular derivative. Independently a similar results have
been obtained by Rodin and Warschawski \cite{Rodin}.

The reduced modulus of a digon is not invariant under conformal mapping.
The following result gives a change-of-variable formula, see, e.g., \cite{VasBOOK}.  Let the digon $D$ with the vertices at $a$ and $b$ be so that the
limit (\ref{eq:moddig}) exists and the Stolz angles are
$\varphi_a$ and $\varphi_b$ . Suppose that there is a conformal
map $f(z)$ of the digon $D$ (which is conformal at $a$, $b$) onto
a digon $D'$, so that there exist the angular limits $f(a)$,
$f(b)$ with the inner angles $\psi_a$ and $\psi_b$ at the vertices
$f(a)$ and $f(b)$ which we also understand as the supremum over
all Stolz angles inscribed in $D'$ with the vertices at $f(a)$ or
$f(b)$ respectively. If the function $f$ has the  finite non-zero
angular derivatives $f'(a)$ and $f'(b)$, then
$\varphi_a=\psi_{a}$, $\varphi_b=\psi_{b}$, and the reduced
modulus (\ref{eq:moddig}) of $D'$ exists and changes {\rm
\cite{Em3}, \cite{Kuz4}, \cite{Sol7}, \cite{VasBOOK}} according to the rule
\begin{equation}
m(f(D),f(a),f(b))= m(D,a,b)+ \frac{1}{\psi_a}\log
|f'(a)|+\frac{1}{\psi_b}\log |f'(b)|.\label{eq:modchange1}
\end{equation}

If we suppose, moreover, that $f$ has the expansion
$$
f(z)=w_1+(z-a)^{\psi_a/\varphi_a}(c_1+c_2(z-a)+\dots)
$$
in a neighborhood of the point $a$, and the expansion
$$
f(z)=w_2+(z-b)^{\psi_b/\varphi_b}(d_1+d_2(z-a)+\dots)
$$
in a neighborhood of the point $b$, then
the reduced modulus of  $D$
changes according to the rule
\begin{equation}
m(f(D),f(a),f(b))= m(D,a,b)+
\frac{1}{\psi_a}\log |c_1|+\frac{1}{\psi_b}\log |d_1|.\label{eq:modchange2}
\end{equation}
Obviously, one can extend this definition to the case of vertices
with infinite support.

\subsection{Extremal partition by digons}

Let $S_0$ be a hyperbolic Riemann surface of finite type.
A finite number $\gamma=(\gamma_1,\dots,
\gamma_m)$
of  simple arcs on $S_0$ that are not freely
homotopic  pairwise on $S_0$ is called an  admissible system
of curves on $S_0$ if these curves are not homotopically trivial, start and finish at
fixed points (or at the same point) which can be either punctures or
else points of hyperbolic components of $S_0$, are not homotopic to  a point
of $S_0$ and do not intersect. A digon $D_j$ on $S_0$ with  two fixed vertices on its
boundary (maybe the same point) is said to be of homotopy type $\gamma_j$
if any arc on $S_0$ connecting two vertices is  homotopic (not
freely) to  $\gamma_j$ within the admissible system.

A system of non-overlapping  digons $(D_1,\dots, D_m)$ on
$S_0$ is said to be of  homotopic type $(\gamma_1,\dots,
\gamma_m)$ if $(\gamma_1,\dots, \gamma_m)$ is an admissible system
of curves on $S_0$ and for any $j\in \{1,\dots,m\}$ the domain
$D_j$ is of homotopic type $\gamma_j$. 

We fix a height vector $\alpha=(\alpha_1,\dots,\alpha_m)$ and require the digons  to be
conformal at their vertices and to satisfy the condition of 
compatibility of angles and heights, e.g., for $a_j$, $b_j$ on hyperbolic smooth components of $S_0$, $\varphi_{a_j}=\pi
\alpha_j/(\sum_{k\in I_{a_j}}\alpha_k)$, $j=1,\dots,m$, $\varphi_{b_j}=\pi
\alpha_j/(\sum_{k\in I_{b_j}}\alpha_k)$, where $I_{a_j}$ ($I_{b_j}$) is the
set of indices which refer to the digons $D_j$ with their vertices
at $a_j$ ($b_j$). With a given admissible system $\gamma$ and a height vector
$\alpha$ we associate the collection $\mathfrak D=(D_1,\dots,
D_m)$ of domains of the homotopy type $\gamma$, satisfying the
condition of compatibility of angles and heights. Such a collection is said to be
associated with $\gamma$ and $\alpha$. Some of digons $(D_1,\dots, D_m)$ (not all of them) may
degenerate. In this case we assume the reduced modulus to be zero.
A general theorem (see \cite{Em3, Kuz4, Sol7, VasBOOK}) implies that
any
collection of non-overlapping admissible digons $\mathfrak D=(D_1,\dots,
D_k)$  associated with the admissible system of
curves $\gamma$ and the vector $\alpha$ satisfies the following
inequality 
\begin{equation}
\sum\limits_{j=1}\sp{m}\alpha_{j}\sp2 m(D_j,a_j,b_j)\geq \sum\limits_{j=1}\sp{m}\alpha_{j}\sp2 m(D_j^*,a_j,b_j), 
\end{equation}
with the equality sign only for $\mathfrak D=\mathfrak D^*$. Here
each $D_j\sp*$ is  a strip domain in the trajectory structure
of a unique quadratic differential $Q(\zeta) d\zeta\sp2$, and
there is a conformal map $g_j(\zeta)$,
$\zeta\in D_j\sp*$ that satisfies the differential equation 
\begin{equation}\label{eq_g}
\alpha_j\sp2 \left( \frac{g'_j(\zeta)}{g_j(\zeta)} \right)\sp2=
4\pi\sp2 Q(\zeta),\quad j=1,\dots,m 
\end{equation} 
 and which maps
$D_j\sp*$ onto the strip $\mathbb C\setminus [0,\infty)$.
The critical trajectories of $Q (\zeta) d\zeta\sp2$ split
$S_0$ into at most $m$  strip domains $\mathfrak
D^*=(D_1^*,\dots,D_m^*)$ associated  with the
admissible system of curves and the height vector (some of
$D_j^*$ can degenerate), satisfying the condition of compatibility of angles and heights.

\section{Extremal problem for one fixed boundary point}

We now let $S_0$ be the Riemann surface $U\setminus\{z\}$, and let $\gamma$ be an admissible curve system
consisting of just one simple arc connecting $z$ with the boundary point 1. We consider a system
of admissible digons $D$ with vertices at $z\in U$ and 1 with inner angles $\varphi_z=2\pi$ and $\varphi_1=\pi$ respectively. Then $D^*_z$ is a strip domain in the trajectory structure of a unique quadratic
differential
$$
Q_z(\zeta)d\zeta^2=A\frac{(\zeta-e^{i\delta})^2}{(\zeta-z)^2(\zeta-1/\bar{z})^2(\zeta-1)^2}d\zeta^2.
$$
The complex constant $A$ is defined by two conditions: the unit circle is a trajectory of $Q_z(\zeta)d\zeta^2$ (the argument)
and the integration of a branch of $\sqrt{Q_z(\zeta)}d\zeta$ along the orthogonal trajectory gives 1 (the absolute value). In this case there is a conformal map $g(\zeta)$,
$\zeta\in D\sp*_z$ that satisfies the differential equation $$
\left( \frac{g'(\zeta)}{g(\zeta)} \right)\sp2=
4\pi\sp2 Q_z(\zeta)$$ and maps
$D\sp*_z$ onto the strip $\mathbb C\setminus [0,\infty)$. To simplify calculation we observe that
the extremal configuration for the above problem is reduced to the problem with $z=0$ by the
M\"obius transformation $$B_z(\zeta)=\frac{1-\bar{z}}{1-z}\frac{\zeta-z}{1-\zeta\bar{z}}.$$ The extremal configuration for $z=0$ is given
by the digon $G=U\setminus (-1,0]$ with its vertices at 0 and 1, and with the reduced modulus $m(G,0,1)=0$.
Making use of the change-of-variable formula (\ref{eq:modchange1}) or (\ref{eq:modchange2}) we find that
$$
m(D^*_z,z,1)=m(G,0,1)+\frac{1}{2\pi}\log |B'_z(0)|+\frac{1}{\pi}\log |B'_z(1)|=\frac{1}{2\pi}\log\frac{(1-|z|^2)^3}{|1-z|^4}.
$$
Let $\phi$ be a univalent conformal mapping $U\to U$ having the angular limit $\phi(1)=1$ and the angular derivative
$\phi(1)=\beta$ fixed in $(0,\infty)$. Then the digon $f(D^*_z)$ is admissible in the problem of the extremal partition
of $S_0$ with changing $z$ to $\phi(z)$ and the domain $D^*_{\phi(z)}$ is the extremal one. So the inequality
$m(\phi(D^*_z),\phi(z),1)\geq m(D^*_{\phi(z)},\phi(z),1)$ holds with the equality sign only for the domain (and the function) $\phi(D^*_z)=D^*_{\phi(z)}$. The inequality is equivalent to
$$
|\phi'(z)|\geq \frac{1}{\beta^2}\frac{(1-|z|^2)^3}{|1-z|^4}\frac{|1-\phi(z)|^4}{(1-|\phi(z)|^2)^3}.
$$

\begin{remark}
We observe that consideration of non-univalent functions gives no result. For example, $\phi(\zeta)=\zeta^2$ satisfies the conditions $\phi(1)=1$, $\phi'(1)=2$, however $\phi'(0)=0$.
\end{remark}

Let us discuss the extremal function. The uniqueness of the extremal configuration in the trajectory structure of the
quadratic differential $Q_z(\zeta)d\zeta^2$ implies that the extremal function $w=\phi^*(\zeta)$ satisfies the differential
equation $Q_z(\zeta)d\zeta^2=Q_{\phi(z)}(w)dw^2$ and the condition $\phi'(1)=\beta$. Since the extremal configuration
is reduced to the problem with $z=0$ by the M\"obius transformation $B_z(\zeta)$, we construct the extremal function
by the following diagram

\begin{figure}[ht]
\centering
\begin{picture}(50,50)
\put(12,5){\vector(1,0){28}} \put(12,45){\vector(1,0){28}}
\put(4,40){\vector(0,-1){30}} \put(49,40){\vector(0,-1){30}}
\put(1,43){$U$} \put(45,43){$U$} \put(45,1){$p_{\alpha}(U)$}\put(-20,0){$\phi^*(U)$}
\put(18,10){$\scriptstyle{B_{\phi^*(z)}}$}\put(20,48){$\scriptstyle{B_z}$}\put(-8,25){$\scriptstyle{\phi^*}$}
\put(37,25){$\scriptstyle{p_{\alpha}}$}
\end{picture}
\end{figure}

\noindent
where $p_{\alpha}(\zeta)$ is  the classical conformal Pick map
$$
p_{\alpha}(z)=\frac{4\alpha z}{\left(1-z+\sqrt{(1-z)^2+4\alpha z}\right)^2}=\alpha z+\dots
$$
of $U$ onto $U\setminus(-1,\,-\alpha/(1+\sqrt{1-\alpha})^2]$. In fact, for a fixed $z\in U$, the number $\alpha$ and the extremal function $\phi^*(z)$ are related by the formula 
$$
\alpha=\frac{1}{\beta^2}\frac{|1-\phi^*(z)|^4(1-|z|^2)^2}{|1-z|^4(1-|\phi^*(z)|^2)^2}.
$$
Thus we proved the following theorem.

\begin{theorem}\label{theoremA}
Let $\phi$ be a conformal univalent map of $U$ into $U$ which is conformal at the boundary point 1 and
$\phi(1)=1$, $\phi'(1)=\beta$. Then,
$$
|\phi'(z)|\geq \frac{1}{\beta^2}\frac{(1-|z|^2)^3}{|1-z|^4}\frac{|1-\phi(z)|^4}{(1-|\phi(z)|^2)^3}.
$$
With a fixed $z\in U$ and $\phi(z)=w$, the equality sign is attained only for the function $\phi^*=B^{-1}_{w}\circ p_{\alpha}\circ B_z$.
\end{theorem}

In view of the results \cite{AndRovnyak, Ruscheweyh} it would be nice to have estimates from below for the $n$-th derivative of a univalent function mapping $U$ into $U$. The present method is closely linked to the first derivative, arising, as it does, from \cite[Theorem 2.2.2, equation (2.8)]{VasBOOK}

\section{Extremal problem for $n$ fixed boundary points}

Let $\gamma$ be an admissible curve system on $S_0=U\setminus \{z\}$ 
consisting of $n$ arcs connecting $0$ with fixed boundary points $\zeta_1,\dots,\zeta_n$. We consider a system
of admissible digons $\mathfrak D=(D_1,\dots, D_n)$ with vertices at $z$ and $\zeta_1,\dots,\zeta_n$ with inner angles $\varphi_{z_j}$ and $\varphi_{\zeta_j}=\pi$ respectively, where $z_j$ is the vertex of $D_j$ supported at $z$. Given a vector $\alpha=(\alpha_1,\dots,\alpha_n)$ normalized by $\sum_{j=1}^n\alpha_j=1$, the condition of compatibility of angles and heights reads as $\varphi_{z_j}=2\pi\alpha_j$. Any collection of non-overlapping
admissible digons $\mathfrak D$ associated with the curve system $\gamma$ and the condition of compatibility of angles and heights satisfies the following inequality
\begin{equation}
\sum\limits_{j=1}^n \alpha_j^2m(D_j,z_j,\zeta_j)\geq \sum\limits_{j=1}^n \alpha_j^2m(D_j^*,z_j,\zeta_j),\label{modulus2}
\end{equation}
for a fixed vector $\alpha$, with the equality sign only for the extremal system of admissible digons $\mathfrak D^*=(D_1^*,\dots,D_n^*)$. We denote this minimum by $\mathcal M(\mathfrak D^*,\gamma,\alpha)$.

Let us consider the case $z=0$. Then the $D^*_j$ are strip domains in the trajectory structure of a unique quadratic
differential
$$
Q(\zeta)d\zeta^2=A\frac{\prod_{k=1}^n(\zeta-e^{i\delta_k})^2}{\zeta^2\prod_{j=1}^n(\zeta-\zeta_j)^2}d\zeta^2.
$$
The complex constant $A$ is positive because of the local trajectory structure about 0, $Q(\zeta)=\frac{A}{\zeta^2}(1+\dots)$. There is a conformal map $g_j(z)$, $z\in D_j^*$ satisfying the differential equation
\begin{equation}\label{eq_g2}
\alpha_j\sp2 \left( \frac{g'_j(\zeta)}{g_j(\zeta)} \right)\sp2=
4\pi\sp2 A\frac{\prod_{k=1}^n(\zeta-e^{i\delta_k})^2}{\zeta^2\prod_{j=1}^n(\zeta-\zeta_j)^2},\quad j=1,\dots,n. 
\end{equation}  
The function $g_j(z)$ maps $D_j^*$ onto the strip $\mathbb C\setminus[0,\infty)$.

The critical trajectories of $Q(\zeta)d\zeta^2$ split $S_0$ into at most $n$ strip domains associated with the admissible system.

The mapping $g_j(\zeta)$ has the expansion $g_j(\zeta)=z^{1/\alpha_j}(c_1+\dots)$ about the origin. Denote $\zeta_j=e^{i\theta_j}$. Letting $\zeta\to 0$
in (\ref{eq_g2}) within $D^*_j$ for any $j=1,\dots, n$, we obtain that $A=1/4\pi^2$ and $\sum_{k=1}^n(\delta_k-\theta_k)$ is of mod($\pi$). Without loss of generality we assume the order $\theta_{k-1}\leq \delta_k\leq \theta_k$, wherehence
$$0<\sum_{k=1}^n(\delta_k-\theta_k)\leq \theta_n-\theta_1<2\pi,$$ therefore it equals $\pi$. The mapping $g_j(\zeta)$ has the expansion $g_j(\zeta)=(z-e^{i\theta_j})^{-2}(d_1+\dots)$ about $\zeta_j$ within $D_j^*$. Letting $\zeta\to \zeta_j$
in (\ref{eq_g2}) within $D^*_j$ for  $j=1,\dots, n$, we obtain a system of equations for $\delta_k$
\begin{equation}
2\alpha_j=\frac{\prod_{k=1}^n(e^{i\theta_j}-e^{i\delta_k})}{e^{i\theta_j}\prod_{k\neq j}^n(e^{i\theta_j}-e^{i\theta_k})},\quad j=1,\dots, n.\label{uu}
\end{equation}

Let $\phi:\,\,U\to U$ is such that $\phi(0)=0$ and $\phi(\zeta_j)=\zeta_j$ for all $j=1,\dots,n$. Then $\phi(\mathfrak D^*)$ is an admissible system of digons
for the same problem of extremal partition satisfying the same condition of compatibility of heights and angles. Therefore,
\begin{equation}\label{uuu1}
\sum\limits_{k=1}^n \alpha_k^2m(\phi(D_k),0_k,\zeta_k)\geq \mathcal M(\mathfrak D^*,\gamma,\alpha).
\end{equation}
The change-of-variable formulas (\ref{eq:modchange1},\ref{eq:modchange2}) imply that
\begin{eqnarray}
& &\sum\limits_{k=1}^n \alpha_k^2 m(\phi(D_k),0_k,\zeta_k)\nonumber \\
&=& \mathcal M(\mathfrak D^*,\gamma,\alpha)+\sum\limits_{k=1}^n \frac{\alpha_k^2}{\pi}\log \phi'(\zeta_k)+\frac{1}{2\pi}\log |\phi'(0)|.\label{uuu2}
\end{eqnarray}
These two relations and the normalization of the vector $\alpha$ yield the inequality
\begin{equation}\label{main}
|\phi'(0)|\geq\frac{1}{\prod_{j=1}^n(\phi'(\zeta_j))^{2\alpha_j^2}},\quad \sum_{j=1}^n\alpha_j=1.
\end{equation}
The extremal configuration is unique and the extremal function $w=\phi^*(\zeta)$ satisfies the complex differential
equation
\begin{equation}
\frac{dw}{d\zeta}=\frac{w\prod_{j=1}^{n}(\zeta-e^{i\delta_j})\prod_{k=1}^{n}(w-\zeta_k)}{\zeta\prod_{j=1}^{n}(w-e^{i\delta_j})\prod_{k=1}^{n}(\zeta-\zeta_k)},\label{eq12}
\end{equation}
normalized as above.
The domain $\phi^*(U)$ is the unit disk $U$ minus at most $n$ analytic arcs starting at the points $e^{i\delta_j}$
along the trajectories of the quadratic differential $Q(w)dw^2$. Their length depends on concrete values of the components of the vector $\alpha$ and on the derivatives $\phi'(\zeta_j)$. We summarize the above in the following theorem.

\begin{theorem}\label{theorem}
Let $\phi$ be a conformal univalent map of $U$ into $U$ which is conformal at the boundary points $\zeta_1\dots,\zeta_n$ and
$\phi(0)=0$, $\phi(\zeta_j)=\zeta_j$, $\phi'(\zeta_j)=\beta_j$, $j=1,\dots, n$. Then for any non-negative vector $\alpha=(\alpha_1,\dots,\alpha_n)$, such that $\sum_{j=1}^n\alpha_j=1$, the following sharp inequality
\[
|\phi'(0)|\geq\frac{1}{\prod_{j=1}^n(\beta_j)^{2\alpha_j^2}},
\]
holds. The equality sign is attained only for the function $\phi^*$ defined by (\ref{eq12}).
\end{theorem}

The minimum in the above inequality with respect to $\alpha$ is attained for
\begin{equation}\label{choice}
\alpha_j=\left(\log \beta_j\sum\limits_{k=1}^n \frac{1}{\log \beta_k}\right)^{-1}.
\end{equation}
This leads to an interesting inequality resembling the Cowen and Pommerenke result \cite{CowenPommerenke} on the derivative at a boundary Denjoy-Wolff point.

\begin{corollary}
Let $\phi$ be a conformal univalent map of $U$ into $U$ which is conformal at the boundary points $\zeta_1\dots,\zeta_n$ and
$\phi(0)=0$, $\phi(\zeta_j)=\zeta_j$, $\phi'(\zeta_j)=\beta_j$, $j=1,\dots, n$. Then the following sharp inequality
\[
\sum_{j=1}^n\frac{1}{\log\beta_j}\leq\frac{-2}{\log|\phi'(0)|}
\]
holds. The equality sign is attained only for the function $\phi^*$ with $\alpha$ chosen as in (\ref{choice}).
\end{corollary}

One may modify this result to the invariant form (wihout fixing 0) by a M\"obius transform, however this requires explicit calculation of moduli as follows. 

Taking into account (\ref{uu}) we see that
\[
\frac{\prod_{k=1}^n(\zeta-e^{i\delta_k})}{\prod_{j=1}^n(\zeta-e^{i\theta_k})}=1+2\sum_{k=1}^n\frac{\alpha_ke^{i\theta_k}}{\zeta-e^{i\theta_k}},
\]
and
\[
\frac{\prod_{k=1}^n(\zeta-e^{i\delta_k})}{\zeta\prod_{j=1}^n(\zeta-e^{i\theta_k})}=-\frac{1}{\zeta}+2\sum_{k=1}^n\frac{\alpha_k}{\zeta-e^{i\theta_k}}.
\]
We take the square root in (\ref{eq_g2}) fixing a branch and integrate the result. Finally, the normalization implies that
\[
g_j(\zeta)=e^{i\varkappa_j}\frac{\zeta^{1/\alpha_j}}{\left(\prod_{k=1}^n(\zeta-e^{i\theta_k})\right)^{1/\alpha_j}}.
\]
The angle $\varkappa_j$ is chosen so that the arc of the unit circle which is the boundary of the domain $D_j^*$ is mapped into the real axis. Obviously, the reduced modulus of the digon $\mathbb C\setminus [0,\infty)$ with respect to its vertices $0$ and $\infty$ vanishes. Applying the change-of-variable formula (\ref{eq:modchange1}) we calculate that
\[
m(D_j^*,0_j,\zeta_j)=\frac{1}{\alpha_j\pi}\log \prod\limits_{k\neq j}|e^{i\theta_j}-e^{i\theta_k}|^{\alpha_k}.
\]

Let us return back to the original problem of the extremal partition (\ref{modulus2}). Applying the M\"obius transformation
\[
M(\zeta)=\frac{\zeta-z}{1-\zeta\bar{z}},
\]
and the change-of-variable formula (\ref{eq:modchange1}), we get
\[
m(D_j^*,z_j,\zeta_j)=\frac{1}{\alpha_j\pi}\log \frac{(1-|z|^2)^{\alpha_j+1/2}}{|1-e^{i\theta_j}\bar{z}|^{2\alpha_j}}\prod\limits_{k\neq j}\bigg|\frac{e^{i\theta_j}-z}{1-e^{i\theta_j}\bar{z}}-\frac{e^{i\theta_k}-z}{1-e^{i\theta_k}\bar{z}}\bigg|^{\alpha_k}.
\]
Analogously to (\ref{uuu1}, \ref{uuu2}) we can prove the following theorem.

\begin{theorem}
Let $\phi$ be a conformal univalent map of $U$ into $U$ which is conformal at the boundary points $\zeta_1\dots,\zeta_n$ and
$\phi(\zeta_j)=\zeta_j$, $\phi'(\zeta_j)=\beta_j$, $j=1,\dots, n$. Then for any non-negative vector $\alpha=(\alpha_1,\dots,\alpha_n)$, such that $\sum_{j=1}^n\alpha_j=1$, the following sharp inequality
\[
|\phi'(z)|\geq\frac{1}{\prod_{j=1}^n(\beta_j)^{2\alpha_j^2}}\prod_{j=1}^n\frac{F_j(z)}{F_j(w)},
\]
holds, where
\[
F_j(z)= \frac{(1-|z|^2)^{\alpha_j(2\alpha_j+1)}}{|1-\zeta_j\bar{z}|^{4\alpha^2_j}}\left(\prod\limits_{k\neq j}\bigg|\frac{\zeta_j-z}{1-\zeta_j\bar{z}}-\frac{\zeta_k-z}{1-\zeta_k\bar{z}}\bigg|^{\alpha_k}\right)^{2\alpha_j}.
\] 
 The equality sign is attained for a function $\phi^*$ constructed analogously to that in Theorem~\ref{theorem}.
\end{theorem}

The details are omitted.

\end{document}